\documentclass[a4paper,11pt]{amsart}
\usepackage{graphicx}
\usepackage{mathptmx}
\usepackage{amsmath}
\usepackage{amssymb}
\usepackage{enumitem}
\usepackage{xcolor}
\usepackage{xparse}
\NewDocumentCommand{\eulerian}{omm}
 {%
  \genfrac<>{0pt}{}{#2}{#3}%
  \IfValueT{#1}{_{\!#1}}%
 }

 \newmuskip\pFqmuskip

\newcommand*\pFq[6][8]{%
  \begingroup 
  \pFqmuskip=#1mu\relax
  \mathchardef\normalcomma=\mathcode`,
  \mathcode`\,=\string"8000
  \begingroup\lccode`\~=`\,
  \lowercase{\endgroup\let~}\pFqcomma
  {}_{#2}F_{#3}{\left(\genfrac..{0pt}{}{#4}{#5}\bigg|#6\right)}%
  \endgroup
}
\newcommand{\pFqcomma}{{\normalcomma}\mskip\pFqmuskip}

\newtheorem{theorem}{Theorem}

\newtheorem{corollary}[theorem]{Corollary}

\begin{document}

\title[A study on properties of degenerate Poisson random variables]{A study on properties of degenerate and zero-truncated degenerate Poisson random variables}

\author{Taekyun  Kim}
\address{Department of Mathematics, Kwangwoon University, Seoul 139-701, Republic of Korea}
\email{tkkim@kw.ac.kr}

\author{DAE SAN KIM}
\address{Department of Mathematics, Sogang University, Seoul 121-742, Republic of Korea}
\email{dskim@sogang.ac.kr}

\author{Hyunseok  Lee}
\address{Department of Mathematics, Kwangwoon University, Seoul 139-701, Republic of Korea}
\email{luciasconstant@kw.ac.kr}

\author{SeongHo Park}
\address{Department of Mathematics, Kwangwoon University, Seoul 139-701, Republic of Korea}
\email{abcd2938471@kw.ac.kr}

\author{Jongkyum Kwon}
\address{Department of Mathematics Education, Gyeongsang National University, Jinju,  Republic of Korea}
\email{mathkjk26@gnu.ac.kr}


\maketitle

\begin{abstract}
In previous papers, studied are the degenerate Poisson random variables and the zero-truncated degenerate Poisson random variables, respectively as degenerate versions of the Poisson random variables and the zero-truncated Poisson random variables. The aim of this paper is to show that various moments of both kinds of random variables can be expressed in terms of certain special polynomials.
\end{abstract}

\section{Introduction}

Carlitz [2] initiated a study on degenerate versions of Bernoulli and Euler numbers which has been extended recently to the researches on various degenerate versions of quite a few special numbers and polynomials.
They have been explored by using several different tools including generating functions, combinatorial methods, $p$-adic analysis, umbral calculus, special functions, differential equations and probability theory as well. \par
The degenerate Poisson random variables are degenerate versions of the Poisson random variables.
In [6], studied are the degenerate binomial and degenerate Poisson random variables in relation to the degenerate Lah-Bell polynomials. Among other things, it is shown that the rising factorial moments of the degenerate Poisson random variable are expressed by the degenerate Lah-Bell polynomials. Also, it is shown that the probability-generating function of the degenerate Poisson random variable is equal to the generating function of the degenerate Lah-Bell polynomials. \par
The zero-truncated Poisson distributions (also called the conditional or the positive Poisson distributions) are certain discrete probability distributions whose supports are the set of positive integers. 
In [10], the zero-truncated degenerate Poisson random variables, whose probability mass functions are a natural extension of the zero-truncated Poisson random variables, are introduced and various properties of those random variables are investigated. Specifically, for those distributions,  studied are its expectation, its variance, its $n$-th moment, its cumulative distribution function and certain expressions for the probability function of a finite sum of independent degenerate zero-truncated Poisson random variables with equal and unequal parameters. \par
As we mentioned in the above, it has been shown that various probabilistic methods can be applied to the study of some special numbers and polynomials arising from combinatorics and number theory. The aim of this paper is to show that various moments for both the degenerate Poisson random variables and the zero-truncated degenerate Poisson random variables can be expressed in terms of various special polynomials. In more detail, for the degenerate Poisson random variables, it is shown that the falling factorial moment is expressed in terms of degenerate Bell polynomials and degenerate Stirling numbers of the first kind, and the $n$th moment is represented by dimorphic degenerate Bell polynomials. As to the zero-truncated degenerate Poisson random variables, it is proved that the rising factorial moment is expressed in terms of the zero-truncated degenerate Lah-Bell polynomials, the falling factorial moment is represented by degenerate Bell polynomials and the degenerate Stirling numbers of the first kind and the $\lambda$ falling factorial moment is given by the degenerate Bell polynomials.\par
For the rest of this section, we recall the necessary facts which are needed throughout this paper.

For any $\lambda\in\mathbb{R}$, the degenerate exponential function is defined as 
\begin{equation}
e_{\lambda}^{x}(t)=\sum_{n=0}^{\infty}(x)_{n,\lambda}\frac{t^{n}}{n!},\quad e_{\lambda}(t)=e_{\lambda}^{1}(t),\quad (\mathrm{see}\ [7,8,9,12]).\label{1}
\end{equation}
where the $\lambda$ falling factorial sequence is given by $(x)_{0,\lambda}=1,\ (x)_{n,\lambda}=x(x-\lambda)\cdots(x-(n-1)\lambda)$, $(n\ge 1)$. \par 
\noindent Let $\log_{\lambda}(t)$ be the compositional inverse of $e_{\lambda}(t)$. Then $\log_{\lambda}(1+t)$ is given by 
\begin{equation}
\log_{\lambda}(1+t)=\sum_{n=1}^{\infty}\lambda^{n-1}(1)_{n,\frac{1}{\lambda}}\frac{t^{n}}{n!}
=\frac{1}{\lambda}\big((1+t)^{\lambda}-1\big),\quad(\mathrm{see}\ [12]).\label{2}	
\end{equation}
The falling and rising factorial sequences are respectively given by 
\begin{align}
&(x)_{0}=1,\ (x)_{n}=x(x-1)\cdots(x-n+1),\ (n\ge 1), \label{3} \\
&\langle x\rangle_{0}=1,\ \langle x \rangle_{n}=x(x+1)\cdots(x+n-1),\ (n\ge 1),
\quad (\mathrm{see}\ [1,2,3,13,14,18]). \nonumber
\end{align} \par
Recently, the degenerate Stirling numbers of the first kind are defined by 
\begin{equation}
(x)_{n}=\sum_{l=0}^{n}S_{1,\lambda}(n,l)(x)_{l,\lambda},\ (n\ge 0), \ \frac{1}{k!}(\log_{\lambda}(1+t))^{k}=\sum_{n=k}^{\infty}S_{1,\lambda}(n,k)\frac{t^{n}}{n!},\ (k\ge 0),\quad (\mathrm{see}\ [7,11]).\label{4}
\end{equation}
Note that $\displaystyle \lim_{\lambda\rightarrow 0}S_{1,\lambda}(n,l)=S_{1}(n,l)\displaystyle$, where $S_{1}(n,l)$ are the ordinary Stirling numbers of the first kind (see [4,16,18]).\par
\noindent As the inversion formula of \eqref{4}, the degenerate Stirling numbers of the second kind are defined as 
\begin{equation}
(x)_{n,\lambda}=\sum_{l=0}^{n}S_{2,\lambda}(n,l)(x)_{l},\ (n\ge 0),\ \frac{1}{k!}(e_{\lambda}(t)-1)^{k}=\sum_{n=k}^{\infty}S_{2,\lambda}(n,k)\frac{t^n}{n!},\ (k \ge 0),\ (\mathrm{see}\ [7,11]). \label{5}
\end{equation} \par
\noindent It is easy to see that the degenerate Stirling numbers of the first and the second kinds satisfy the following orthogonality relations:
\begin{equation}
\sum_{l=0}^{n}S_{1,\lambda}(n,l)S_{2,\lambda}(l,k)=\sum_{l=0}^{n}S_{2,\lambda}(n,l)S_{1,\lambda}(l,k)=\delta_{n,k}. \label{6}
\end{equation}
It is known that $X=X_{\lambda}$ is the degenerate Poisson random variable with parameter $\alpha>0$, which is denoted by $X \sim \mathrm{Poi}_{\lambda}(\alpha)$, if the probability mass function of $X$ is given by 
\begin{equation}
	p(i)=P\{X=i\}=e_{\lambda}^{-1}(\alpha)\cdot\frac{\alpha^{i}}{i!}(1)_{i,\lambda},\quad (i\ge 0),\quad (\mathrm{see}\ [10,11]). \label{7}
\end{equation}
Let $f$ be a real valued function. Then the expectation of $f(X)$ is defined by 
\begin{equation}
E\big[f(X)\big]=\sum_{n=0}^{\infty}f(n)p(n)=\sum_{n=0}^{\infty}f(x)P\{X=n\},\quad (\mathrm{see}\ [15,17]),\label{8}	
\end{equation}
where $p(n)$ is the probability mass function of $X$. \par 
Let $X$ be the degenerate Poisson random variable with parameter $\alpha(>0)$. Then the generating function of the moment of $X$ is given by 
\begin{align}
E\Big[e_{\lambda}^{X}(t)\Big]\ &=\ \sum_{i=0}^{\infty}e_{\lambda}^{i}(t)\cdot p(i)\ =\ e_{\lambda}^{-1}(\alpha)\sum_{i=0}^{\infty}e_{\lambda}^{i}(t)\frac{\alpha^{i}}{i!}(1)_{i,\lambda} \label{9}\\
&=\ e_{\lambda}^{-1}(\alpha)e_{\lambda}\Big(\alpha e_{\lambda}(t)\Big)\ =\ \sum_{n=0}^{\infty}\mathrm{Bel}_{n,\lambda}(\alpha)\frac{t^{n}}{n!},\quad (\mathrm{see}\ [6,10]),\nonumber	
\end{align}
where $\mathrm{Bel}_{n,\lambda}(x)$ are the degenerate Bell polynomials given by
\begin{equation}
e_{\lambda}^{-1}(x)e_{\lambda}\Big(x e_{\lambda}(t)\Big)\ =\ \sum_{n=0}^{\infty}\mathrm{Bel}_{n,\lambda}(x)\frac{t^{n}}{n!},\quad (\mathrm{see}\ [6,8]).\label{10}
\end{equation}
Note that 
\begin{align}
\sum_{n=0}^{\infty}\lim_{\lambda\rightarrow 0}\mathrm{Bel}_{n,\lambda}(x)\frac{t^{n}}{n!}\ &=\ \lim_{\lambda\rightarrow 0}e_{\lambda}^{-1}(x)e_{\lambda}\Big(xe_{\lambda}(t)\Big)\ =\ e^{x(e^{t}-1)} \label{11} \\
&=\ \sum_{n=0}^{\infty}\mathrm{Bel}_{n}(x)\frac{t^{n}}{n!},\nonumber	
\end{align}
where $\mathrm{Bel}_{n}(x)$ are the ordinary Bell polynomials. \par 
\noindent From \eqref{10}, we note that 
\begin{displaymath}
\mathrm{Bel}_{n,\lambda}(x)=\frac{1}{e_{\lambda}(x)}\sum_{k=0}^{\infty}\frac{x^{k}}{k!}(1)_{k,\lambda}(k)_{n,\lambda},\quad (n\ge 0).
\end{displaymath} 
By \eqref{9}, we get 
\begin{equation}
E[(X)_{n,\lambda}]=\mathrm{Bel}_{n,\lambda}(\alpha), \quad (n\ge 0). \label{12}
\end{equation} \par
It is well known that the Lah-numbers are defined by 
\begin{equation}
\frac{1}{k!}\bigg(\frac{t}{1-t}\bigg)^{k}=\sum_{n=k}^{\infty}L(n,k)\frac{t^{n}}{n!},\quad (k\ge 0),\quad (\mathrm{see}\ [5]). \label{13}	
\end{equation}
From \eqref{13}, the Lah-Bell polynomials are given by the generating function to be 
\begin{equation}
e^{x(\frac{t}{1-t})}=e^{x(\frac{1}{1-t}-1)}=\sum_{n=0}^{\infty}B_{n}^{L}(x)\frac{t^{n}}{n!}. \label{14}
\end{equation}
Note that 
\begin{equation}
B_{n}^{L}(x)=\sum_{k=0}^{n}L(n,k)x^{k},\quad (n\ge 0),\quad (\mathrm{see}\ [5]). \label{15}
\end{equation}
Recently, the degenerate Lah-Bell polynomials are defined as
\begin{equation}
e_{\lambda}^{-1}(x)\cdot e_{\lambda}\bigg(\frac{x}{1-t}\bigg)=\sum_{n=0}^{\infty}B_{n,\lambda}^{L}(x)\frac{t^{n}}{n!},\quad(\mathrm{see}\ [6]). \label{16}
\end{equation}
Note that $\displaystyle\lim_{\lambda\rightarrow 0}B_{n,\lambda}^{L}(x)=B_{n}^{L}(x),\ (n\ge 0)\displaystyle$. \par

\section{Properties of degenerate Poisson Random Variables}
First, we consider the fully degenerate Bell polynomials which are given by 
\begin{equation}
e_{\lambda}^{x}\big(e_{\lambda}(t)-1\big)=\sum_{n=0}^{\infty}\beta_{n,\lambda}(x)\frac{t^{n}}{n!}. \label{17}
\end{equation}
When $x=1$, $\beta_{n,\lambda}=\beta_{n,\lambda}(1)$ are called the fully degenerate Bell numbers. \par 
Replacing $t$ by $\log_{\lambda}(1+t)$ in \eqref{17}, we get 
\begin{align}
e_{\lambda}^{x}(t)\ &= \ \sum_{k=0}^{\infty}\beta_{k,\lambda}(x)\frac{1}{k!}\big(\log_{\lambda}(1+t)\big)^{k}\ =\ \sum_{k=0}^{\infty}\beta_{k,\lambda}(x)\sum_{n=k}^{\infty}S_{1,\lambda}(n,k)\frac{t^{n}}{n!}\label{18} \\
&=\ \sum_{n=0}^{\infty}\bigg(\sum_{k=0}^{n}\beta_{k,\lambda}(x)S_{1,\lambda}(n,k)\bigg)\frac{t^{n}}{n!}.\nonumber
\end{align}
Thus, by \eqref{18} and \eqref{1}, we obtain the following theorem. 
\begin{theorem}
For $n\ge 0$, we have 
\begin{displaymath}
(x)_{n,\lambda}=\sum_{k=0}^{n}\beta_{k,\lambda}(x)S_{1,\lambda}(n,k). 
\end{displaymath}
\end{theorem}
	
By Theorem 1 and \eqref{6}, we obtain the following corollary. 
\begin{corollary}
	For $n\ge 0$, we have 
	\begin{displaymath}
	\beta_{n,\lambda}(x)=\sum_{k=0}^{n}S_{2,\lambda}(n,k)(x)_{k,\lambda}.
	\end{displaymath}
\end{corollary}
Let $X=X_{\lambda}$ be the degenerate Poisson random variable with parameter $\alpha(>0)$. Then we have 
\begin{align}
E\big[(1+t)^{X}\big]\ &=\ E\big[e_{\lambda}^{X}(\log_{\lambda}(1+t))\big] \nonumber \\
&=\ \sum_{k=0}^{\infty}E\big[(X)_{k,\lambda}\big]\frac{1}{k!}\big(\log_{\lambda}(1+t)\big)^{k} \nonumber \\
&=\ \sum_{k=0}^{\infty}E\big[(X)_{k,\lambda}\big]\sum_{n=k}^{\infty}S_{1,\lambda}(n,k)\frac{t^{n}}{n!} \label{19} \\
&=\ \sum_{n=0}^{\infty}\bigg(\sum_{k=0}^{n}E[(X)_{k,\lambda}]S_{1,\lambda}(n,k)\bigg)\frac{t^{n}}{n!}. \nonumber	
\end{align}
Thus, by \eqref{19}, we obtain the following theorem. 
\begin{theorem}
	For $n\ge 0$, we have 
	\begin{displaymath}
		E\big[(X)_{n}\big]=\sum_{k=0}^{n}E\big[(X)_{k,\lambda}\big]S_{1,\lambda}(n,k)=\sum_{k=0}^{n}\mathrm{Bel}_{k,\lambda}(\alpha)S_{1,\lambda}(n,k), 
	\end{displaymath}
	where $X$ is the degenerate Poisson random variable with parameter $\alpha>0$. 
\end{theorem}
For $X\sim\mathrm{Poi}_{\lambda}(\alpha)$, we have 
\begin{align}
E\big[(1+t)^{X}\big]\ &=\ \sum_{i=0}^{\infty}(1+t)^{i}p(i)\ =\ \sum_{i=0}^{\infty}\sum_{m=0}^{i}\binom{i}{m}p(i)t^{m}\label{20}\\
&=\ \sum_{m=0}^{\infty}\bigg(\sum_{i=m}^{\infty}\binom{i}{m}p(i)t^{m}\bigg)\ =\ \sum_{m=0}^{\infty}\frac{e_{\lambda}^{-1}(\alpha)}{m!}\sum_{i=m}^{\infty}\frac{\alpha^{i}(1)_{i,\lambda}}{(i-m)!}t^{m} \nonumber \\
&=\ \sum_{m=0}^{\infty}\frac{e_{\lambda}^{-1}(\alpha)}{m!}\alpha^{m}\bigg(\sum_{i=0}^{\infty}\frac{\alpha^{i}}{i!}(1)_{i+m,\lambda}\bigg)t^{m} \nonumber \\
&=\ \sum_{m=0}^{\infty}\frac{\alpha^{m}}{m!}e_{\lambda}^{-1}(\alpha)(1)_{m,\lambda}\bigg(\sum_{i=0}^{\infty}\frac{(1-\lambda m)_{i,\lambda}}{i!}\alpha^{i}\bigg)t^{m}\nonumber \\
&=\ \sum_{m=0}^{\infty}\bigg(\alpha^{m}e_{\lambda}^{-1}(\alpha)(1)_{m,\lambda}e_{\lambda}^{1-\lambda m}(\alpha)\bigg)\frac{t^{m}}{m!} \nonumber \\
&=\ \sum_{m=0}^{\infty}\bigg(\alpha^{m}(1)_{m,\lambda}\frac{1}{(1+\lambda\alpha)^{m}}\bigg)\frac{t^{m}}{m!}. \nonumber
\end{align}
Therefore, by \eqref{20} and Theorem 3 we obtain the following theorem. 
\begin{theorem}
	For $X\sim\mathrm{Poi}_{\lambda}(\alpha)$, and $m\ge 0$, we have 
	\begin{displaymath}
		E\big[(X)_{m}\big]=\alpha^{m}(1)_{m,\lambda}\frac{1}{(1+\lambda\alpha)^{m}}. 
	\end{displaymath}
	In particular, 
	\begin{displaymath}
		(1+\lambda\alpha)^{m}\sum_{k=0}^{m}\mathrm{Bel}_{k,\lambda}(\alpha)S_{1,\lambda}(m,k)=\alpha^{m}(1)_{m,\lambda}. 
	\end{displaymath}
\end{theorem}
For $X\sim\mathrm{Poi}_{\lambda}(\alpha)$, we observe that 
\begin{align}
	E\bigg[\bigg(\frac{1}{1-t}\bigg)^{X}\bigg]\ &=\ \sum_{i=0}^{\infty}\bigg(\frac{1}{1-t}\bigg)^{i}p(i)\ =\ \sum_{i=0}^{\infty}\bigg(\frac{1}{1-t}\bigg)^{i}e_{\lambda}^{-1}(\alpha)\frac{\alpha^{i}}{i!}(1)_{i,\lambda}\label{21}\\
	&=\ e_{\lambda}^{-1}(\alpha)e_{\lambda}\bigg(\frac{\alpha}{1-t}\bigg)\ =\ \sum_{n=0}^{\infty}B_{n,\lambda}^{L}(\alpha)\frac{t^{n}}{n!}. \nonumber 
\end{align}
On the other hand 
\begin{equation}
E\bigg[\bigg(\frac{1}{1-t}\bigg)^{X}\bigg]=\sum_{n=0}^{\infty}E\bigg[\binom{X+n-1}{n}\bigg]t^{n}\label{22}.
\end{equation}
Therefore, by \eqref{21} and \eqref{22}, we obtain the following theorem. 
\begin{theorem}
	For $X\sim\mathrm{Poi}_{\lambda}(\alpha)$, and $n\ge 0$, we have 
	\begin{displaymath}
		E\bigg[\binom{X+n-1}{n}\bigg]=\frac{1}{n!}B_{n,\lambda}^{L}(\alpha). 
	\end{displaymath}
\end{theorem}
Now, we observe that 
\begin{equation}
\begin{aligned}
E\Big[e^{Xt}\Big]\ &=\ \sum_{i=0}^{\infty}e^{it}p(i)\ =\ e_{\lambda}^{-1}(\alpha)\sum_{i=0}^{\infty}e^{it}\frac{\alpha^{i}}{i!}(1)_{i,\lambda} \\
&=\ e^{-1}_{\lambda}(\alpha)e_{\lambda}(\alpha e^{t})\ =\ \sum_{n=0}^{\infty}\bigg(\frac{1}{e_{\lambda}(\alpha)}\sum_{k=0}^{\infty}\frac{(1)_{k,\lambda}}{k!}\alpha^{k}k^{n}\bigg)\frac{t^{n}	}{n!}, 
\end{aligned}\label{23}
\end{equation}
where $X\sim\mathrm{Poi}_{\lambda}(\alpha)$. \par 
Thus, by \eqref{22}, we get 
\begin{equation}
E[X^{n}]=\frac{1}{e_{\lambda}(\alpha)}\sum_{k=0}^{\infty}\frac{(1)_{k,\lambda}}{k!}\alpha^{k}k^{n},\quad (n\ge 0), \label{24}
\end{equation}
where $X\sim\mathrm{Poi}_{\lambda}(\alpha)$.\par 
Let us define the {\it{dimorphic degenerate Bell polynomials}} as follows: 
\begin{equation}
B_{n,\lambda}(x)=\frac{1}{e_{\lambda}(x)}\sum_{k=0}^{\infty}\frac{(1)_{k,\lambda}}{k!}x^{k}k^{n},\quad (n\ge 0). \label{25}	
\end{equation}
Then, for $X\sim\mathrm{Poi}_{\lambda}(\alpha)$, we get  
\begin{equation}
E[X^{n}]=B_{n,\lambda}(\alpha),\quad (n\ge 0). \label{26}
\end{equation}
Also, for $X\sim\mathrm{Poi}_{\lambda}(\alpha)$, we have
\begin{align}
E\bigg[\bigg(\frac{1}{1-t}\bigg)^{X}\bigg]\ &=\ \sum_{k=0}^{\infty}E[X^{k}](-1)^{k}\frac{1}{k!}\big(\log(1-t)\big)^{k} \label{27} \\
&=\ \sum_{k=0}^{\infty}B_{k,\lambda}(\alpha)(-1)^{k}\sum_{n=k}^{\infty}(-1)^{n}S_{1}(n,k)\frac{t^{n}}{n!}.\nonumber \\
&=\ \sum_{n=0}^{\infty}\bigg(\sum_{k=0}^{n}B_{k,\lambda}(\alpha)(-1)^{n-k}S_{1}(n,k)\bigg)\frac{t^{n}}{n!}. \nonumber 
\end{align}
Therefore, by Theorem 5 and \eqref{27}, we obtain the following theorem. 
\begin{theorem}
	For $X\sim\mathrm{Poi}_{\lambda}(\alpha)$, we have
	\begin{displaymath}
		E[\langle X\rangle_{n}]=\sum_{k=0}^{n}B_{k,\lambda}(\alpha)\big|S_{1}(n,k)\big|,
	\end{displaymath}
	where $|S_{1}(n,k)|$ are the unsigned Stirling numbers of the first kind. \par 
	In particular, we have
	\begin{displaymath}
	B_{n,\lambda}^{L}(\alpha)=\sum_{k=0}^{n}B_{k,\lambda}(\alpha)\big|S_{1}(n,k)\big|.
	\end{displaymath}
\end{theorem}

\section{Properties of zero-truncated degenerate Poisson random variables} 
In this section, we consider the zero-truncated degenerate Poisson random variables with parameter $\alpha>0$. A random variable $X=X_{\lambda}$ is the zero-truncated degenerate Poisson random variables with parameter $\alpha>0$, which is denoted by $X\sim \mathrm{Poi}_{\lambda}^{*}(\alpha)$, if the probability mass function of $X$ is given by 
\begin{equation}
p(k)=P\{X=k\}=\frac{1}{e_{\lambda}(\alpha)-1}\frac{\alpha^{k}}{k!}(1)_{k,\lambda},\quad (\mathrm{see}\ [10]),\label{28}
\end{equation}
where $k=1,2,3,\dots$.\par 
For $X\sim \mathrm{Poi}_{\lambda}^{*}(\alpha)$, by \eqref{28}, we get 
\begin{align}
E\bigg[\bigg(\frac{1}{1-t}\bigg)^{X}\bigg]\ &=\ \sum_{i=1}^{\infty}\bigg(\frac{1}{1-t}\bigg)^{i}p(i)\ =\ \sum_{i=1}^{\infty}\bigg(\frac{1}{1-t}\bigg)^{i}\frac{1}{e_{\lambda}(\alpha)-1}\frac{\alpha^{i}}{i!}(1)_{i,\lambda} \label{29} \\
&=\ \frac{1}{e_{\lambda}(\alpha)-1}\bigg(\sum_{i=0}^{\infty}\bigg(\frac{1}{1-t}\bigg)^{i}\frac{\alpha^{i}}{i!}(1)_{i,\lambda}-1\bigg) \nonumber \\
&=\ \frac{1}{e_{\lambda}(\alpha)-1}\bigg(e_{\lambda}\bigg(\frac{\alpha}{1-t}\bigg)-1\bigg). \nonumber
\end{align}
Now, we consider the zero-truncated degenerate Lah-Bell polynomials given by  
\begin{equation}
\frac{1}{e_{\lambda}(x)-1}\bigg(e_{\lambda}\bigg(\frac{x}{1-t}\bigg)-1\bigg)=\sum_{n=0}^{\infty}B_{n,\lambda}^{(L,0)}(x)\frac{t^{n}}{n!}. \label{30}
\end{equation}
From \eqref{30}, we note that 
\begin{align}
\frac{1}{e_{\lambda}(x)-1}\bigg(e_{\lambda}\bigg(\frac{x}{1-t}\bigg)-1\bigg)\ &=\ \frac{1}{1-e_{\lambda}^{-1}(x)}\bigg(e_{\lambda}^{-1}(x)e_{\lambda}\bigg(\frac{x}{1-t}\bigg)-e_{\lambda}^{-1}(x)\bigg) \label{31}\\ 
&=\ \sum_{n=1}^{\infty}\frac{1}{1-e_{\lambda}^{-1}(x)}B_{n,\lambda}^{L}(x)\frac{t^{n}}{n!}+1.\nonumber 
\end{align}

Therefore, by \eqref{30} and \eqref{31}, we obtain the following theorem. 
\begin{theorem}
For $n\in\mathbb{N}$, we have 
\begin{displaymath}
B_{n,\lambda}^{(L,0)}(x)=\frac{1}{1-e_{\lambda}^{-1}(x)}B_{n,\lambda}^{L}(x),\quad and \quad B_{0,\lambda}^{(L,0)}(x)=1. 
\end{displaymath}
\end{theorem}
From \eqref{29} we obtain the following corollary. 
\begin{corollary}
For $X\sim \mathrm{Poi}_{\lambda}^{*}(\alpha)$, and $n\in\mathbb{N}$, we have 
\begin{displaymath}
E\bigg[\binom{X+n-1}{n}\bigg]=\frac{1}{n!}B_{n,\lambda}^{(L,0)}(\alpha). 
\end{displaymath}
\end{corollary}
For $X\sim \mathrm{Poi}_{\lambda}^{*}(\alpha)$, we observe that 
\begin{align}
\sum_{n=0}^{\infty}E[(X)_{n}]\frac{t^{n}}{n!}\ &=\ E[(1+t)^{X}]\ =\ \frac{1}{e_{\lambda}(\alpha)-1}\sum_{j=1}^{\infty}(1+t)^{j}\frac{\alpha^{j}}{j!}(1)_{j,\lambda}\label{32} \\
&=\ 	\frac{1}{e_{\lambda}(\alpha)-1}\sum_{j=0}^{\infty}(1+t)^{j}\frac{\alpha^{j}}{j!}(1)_{j,\lambda}-\frac{1}{e_{\lambda}(\alpha)-1}\nonumber \\
&=\ \frac{1}{e_{\lambda}(\alpha)-1}\sum_{n=0}^{\infty}\sum_{j=n}^{\infty}\binom{j}{n}\frac{\alpha^{j}}{j!}(1)_{j,\lambda}t^{n}-\frac{1}{e_{\lambda}(\alpha)-1}\nonumber \\
&=\ \frac{1}{e_{\lambda}(\alpha)-1}\sum_{n=0}^{\infty}\frac{1}{n!}\sum_{j=n}^{\infty}\frac{\alpha^{j}}{(j-n)!}(1)_{j,\lambda}t^{n}-\frac{1}{e_{\lambda}(\alpha)-1}\nonumber \\
&=\ \sum_{n=0}^{\infty}\frac{\alpha^{n}}{n!(e_{\lambda}(\alpha)-1)}\sum_{j=0}^{\infty}\frac{\alpha^{j}}{j!}(1)_{j+n,\lambda}t^{n}-\frac{1}{e_{\lambda}(\alpha)-1}\nonumber \\
&=\ \sum_{n=1}^{\infty}\frac{\alpha^{n}(1)_{n,\lambda}}{n!(e_{\lambda}(\alpha)-1)}\sum_{j=0}^{\infty}\frac{\alpha^{j}}{j!}(1-n\lambda)_{j,\lambda}t^{n}+1\nonumber \\
&=\ \sum_{n=1}^{\infty}\frac{\alpha^{n}e_{\lambda}(\alpha)}{n!(e_{\lambda}(\alpha)-1)}(1)_{n,\lambda}e_{\lambda}^{-n\lambda}(\alpha)t^{n}+1\nonumber \\
&=\ \sum_{n=1}^{\infty}\frac{e_{\lambda}^{-n\lambda}(\alpha)\alpha^{n}(1)_{n,\lambda}}{n!(1-e_{\lambda}^{-1}(\alpha)}t^{n}+1.\nonumber 
\end{align}
Therefore, by comparing the coefficients on both sides of \eqref{32}, we obtain the following theorem. 
\begin{theorem}
	For $n\in\mathbb{N}$, $X\sim \mathrm{Poi}^{*}(\alpha)$, we have 
	\begin{displaymath}
		E[(X)_{n}]=\frac{\alpha^{n}(1)_{n,\lambda}}{(1-e_{\lambda}^{-1}(\alpha))}\frac{1}{(1+\lambda\alpha)^{n}}.
	\end{displaymath}
\end{theorem}
For $X\sim \mathrm{Poi}^{*}(\alpha)$, we have 
\begin{align}
&\sum_{n=0}^{\infty}E\big[(X)_{n,\lambda}\big]\frac{t^{n}}{n!} = E\big[e_{\lambda}^{X}(t)\big] = \frac{1}{e_{\lambda}(\alpha)-1}\sum_{i=1}^{\infty}e_{\lambda}^{i}(t)\frac{\alpha^{i}(1)_{i,\lambda}}{i!}\label{33}\\
& = \frac{1}{e_{\lambda}(\alpha)-1}\sum_{i=0}^{\infty}e_{\lambda}^{i}(t)\frac{\alpha^{i}(1)_{i,\lambda}}{i!}-\frac{1}{e_{\lambda}(\alpha)-1} = \frac{1}{e_{\lambda}(\alpha)-1}e_{\lambda}\big(\alpha e_{\lambda}(t)\big)-\frac{1}{e_{\lambda}(\alpha)-1}\nonumber\\
& = \frac{1}{1-e_{\lambda}^{-1}(\alpha)}\cdot e_{\lambda}^{-1}(\alpha)e_{\lambda}\big(\alpha e_{\lambda}(t)\big)-\frac{1}{e_{\lambda}(\alpha)-1} = \frac{1}{1-e_{\lambda}^{-1}(\alpha)}\sum_{n=0}^{\infty}\mathrm{Bel}_{n,\lambda}(\alpha)\frac{t^{n}}{n!}-\frac{1}{e_{\lambda}(\alpha)-1}\nonumber \\
&=\sum_{n=1}^{\infty}\frac{1}{1-e_{\lambda}^{-1}(\alpha)}\mathrm{Bel}_{n,\lambda}(\alpha)\frac{t^{n}}{n!}+1. \nonumber 
\end{align}
Therefore, by comparing the coefficients on both sides of \eqref{33}, we obtain the following theorem. 
\begin{theorem}
For $X\sim \mathrm{Poi}^{*}(\alpha)$, and $n\in\mathbb{N}$, we have 
\begin{displaymath}
E\big[(X)_{n,\lambda}\big]=\frac{1}{1-e_{\lambda}^{-1}(\alpha)}\mathrm{Bel}_{n,\lambda}(\alpha).
\end{displaymath}
\end{theorem}
For $X\sim \mathrm{Poi}^{*}(\alpha)$, we have 
\begin{align}
E\big[(1+t)^{X}\big]\ &=\ E\big[e_{\lambda}^{X}(\log_{\lambda}(1+t))\big]\ =\ \sum_{k=0}^{\infty}E\big[(X)_{k,\lambda}\big]\frac{1}{k!}\big(\log_{\lambda}(1+t)\big)^{k}\label{34} \\
&=\ \sum_{n=0}^{\infty}\bigg(\sum_{k=0}^{n}E\big[(X)_{k,\lambda}\big]S_{1,\lambda}(n,k)\bigg)	\frac{t^{n}}{n!}\nonumber \\
&=\ \sum_{n=1}^{\infty}\bigg(\sum_{k=1}^{n}\frac{1}{1-e_{\lambda}^{-1}(\alpha)}\mathrm{Bel}_{k,\lambda}(\alpha)S_{1,\lambda}(n,k)\bigg)\frac{t^{n}}{n!}+1.\nonumber 
\end{align}
Thus, by \eqref{34}, we get the following result.
\begin{theorem}
For $X\sim \mathrm{Poi}^{*}(\alpha)$, and $n\in\mathbb{N}$, we have 
\begin{displaymath}
E\big[(X)_{n}\big]\ =\ \frac{1}{1-e_{\lambda}^{-1}(\alpha)}\sum_{k=1}^{n}\mathrm{Bel}_{k,\lambda}(\alpha)S_{1,\lambda}(n,k). 
\end{displaymath}
\end{theorem}

\section{Conclusion} 

In recent years, we have seen that degenerate versions of many special polynomials and numbers can be
investigated by means of various different tools. Here we applied probabilistic methods in order to study both the degenerate Poisson random variables and the zero-truncated degenerate Poisson random variables. In more detail, for the degenerate Poisson random variables, we showed that the falling factorial moment is expressed in terms of degenerate Bell polynomials and degenerate Stirling numbers of the first kind, and the $n$th moment is represented by dimorphic degenerate Bell polynomials. For the zero-truncated degenerate Poisson random variables, we proved that the rising factorial moment is expressed in terms of the zero-truncated degenerate Lah-Bell polynomials, the falling factorial moment is represented by degenerate Bell polynomials and the degenerate Stirling numbers of the first kind and the $\lambda$ falling factorial moment is given by the degenerate Bell polynomials.\par
It is one of our future projects to continue this line of research, namely to explore applications of various methods of probability theory to the study of some special polynomials and numbers.

\vspace{0.2in}

\noindent{\bf{Acknowledgments:}} 
Not applicable.
\vspace{0.2in}

\noindent{\bf{Funding:}} 
Not applicable.

\vspace{0.2in}

\noindent{\bf {Availability of data and materials:}}
Not applicable.

\vspace{0.1in}

\noindent{\bf {Competing interests:}}
The authors declare no conflict of interest.

\vspace{0.1in}

\noindent{\bf{Authors' contributions:}} T.K. and D.S.K. conceived of the framework and structured the whole paper; T. K. and D.S.K. wrote the paper; D.S.K. and T.K. completed the revision of the article; H.L, J.K. and
S.P checked the errors of the article.. All authors have read and agreed to the published version of the manuscript.

\vspace{0.1in}


\end{document}